\input amstex 
\documentstyle{amsppt}
\input bull-ppt
\keyedby{bull318/lbr}

\topmatter
\cvol{27}
\cvolyear{1992}
\cmonth{October}
\cyear{1992}
\cvolno{2}
\cpgs{279-283}
\title Keller's cube-tiling conjecture\\is false in high 
dimensions \endtitle
\author Jeffrey C. Lagarias and Peter W. Shor\endauthor
\shortauthor{J. C. Lagarias and P. W. Shor}
\shorttitle{Keller's cube-tiling conjecture}
\address AT\&T Bell Laboratories, Murray Hill, New Jersey 
07974\endaddress
\date March 5, 1992\enddate
\subjclass Primary 52C22 Secondary 05B45\endsubjclass
\abstract O. H. Keller conjectured in 1930 that in any 
tiling of $\Bbb R^n$ by
unit $n$-cubes there exist two of them having a complete 
facet in common. O.
Perron proved this conjecture for $n\le 6$. We show that 
for all $n\ge 10$
there exists a tiling of $\Bbb R^n$ by unit $n$-cubes such 
that no two
$n$-cubes have a complete facet in common.\endabstract
\endtopmatter

\document
\heading 1. Introduction\endheading
In 1907 Minkowski  \cite5 conjectured that all the 
extremal lattices for the
supremum norm were of a certain simple form and observed 
that this conjecture
had a geometric interpretation: in any lattice tiling of 
$\Bbb R^n$ with unit
$n$-cubes there must exist two cubes having a complete 
facet $((n-1)$-face) in
common. He proved this for $n=2$ and $3$. In studying this 
question, Keller
 \cite4 generalized it to conjecture that any tiling of 
$\Bbb R^n$ by
unit $n$-cubes contains two cubes having a complete facet 
in common. In 1940
Perron
 \cite6 proved Keller's conjecture for dimensions $n\le 
6$. Soon after,
Haj\'os  \cite2 proved that Minkowski's original 
conjecture is true in
all dimensions. Keller's stronger conjecture remained 
open. Haj\'os 
\cite3 later gave a combinatorial problem concerning 
factorization of abelian groups,
which he proved was equivalent to Keller's conjecture. 
Stein  \cite7 gave
a survey of these results and other related tiling 
problems. 

More recently, Szab\'o  \cite8 showed that if Keller's 
conjecture is
false in $\Bbb R^n$, then there exists a counterexample 
tiling in some $\Bbb
R^m$ (with $m\ge n)$ having the following extra 
properties: the centers of all
cubes are in $\tfrac12 \Bbb Z^m$, and the tiling is 
periodic with period
lattice containing $2\Bbb Z^m$. Corr\'adi and Szab\'o  
\cite1 studied a
graph-theoretic version of this latter problem, showing 
directly that there are
no such counterexamples for $m\le 5$. 

We explicitly construct a counterexample tiling of 
Szab\'o's type in $\Bbb
R^{10}$.  Keller's conjecture is then false for all $n\ge 
10$ because a
counterexample tiling in $\Bbb R^n$ gives one in $\Bbb 
R^{n+1}$ by ``stacking''
layers of this tiling with suitable translations made 
between adjacent layers.

\heading 2. Main result\endheading

We prove the following result. 

\thm{Theorem A} For $n=10$ and $12$ there exists a tiling 
of $\Bbb R^n$ by unit
cubes such that 
\roster
\item The centers of all cubes are in $\tfrac12\Bbb Z^n$\RM;
\item The tiling is periodic with period lattice $2\Bbb 
Z^n$\RM;
\item No two cubes have a complete facet in common\/.
\endroster
\ethm

Before giving the constructions, we describe Corr\'adi and 
Szab\'o's equivalent
graph-theoretic criterion for such a tiling to exist in 
$\Bbb R^n$. 

Scale everything up by a factor of $2$ to consider tilings 
of $\Bbb R^n$ by
translates of the cube 
$$
\sans C=\{(x_1,\dotsc, x_n):-1\le x_i\le 1\text{ for all 
}i\}
$$
of side 2 centered at the origin, such that the centers of 
all cubes are in $
\Bbb Z^n$, and the tiling is periodic with period lattice 
$4\Bbb Z^n$. There
are $2^n$ equivalence classes of cubes $\bold m+\sans C+
4\Bbb Z^n$ in such a
tiling, and each equivalence class contains a unique cube 
with center 
$$
\bold m=(m_1,\dotsc, m_n)\in \Bbb Z^n,\qquad 0\le m_i\le 
3\.\tag1
$$
The collection $\scr S$ of these $2^n$ vectors describes 
the tiling. 

Now form two graphs $G_n$ and $G^*_n$, each of which has 
$4^n$ vertices
labeled by the $4^n$ vectors in $\Bbb Z^n$ of form 
\thetag1, as follows.
Consider the conditions:
\roster
\item"(a)" \<$\bold m$ and $\bold m'$ have some 
$|m_i-m'_i|=2$. 
\item"(b)" \<$\bold m$ and $\bold m'$ differ in two 
coordinate directions. 
\endroster
$G_n$ has an edge between vertices $\bold m$ and $\bold 
m'$ if (a) holds, while
$G^*_n$ has an edge between $\bold m$ and $\bold m'$ if 
(a) and (b) both hold.
Condition (a) says that all translates under $4\Bbb Z^n$ 
of cubes $\sans C$
centered at $\bold m$ and $\bold m'$ have disjoint 
interiors, while (a) and (b)
together say that all translates under $4\Bbb Z^n$ of such 
cubes also do not
have a complete facet in common. 

A set $\scr S$ of $2^n$ vectors satisfying \thetag1 yields 
a $4\Bbb
Z^n$-periodic cube tiling if and only if $\scr S$ forms a 
clique in $G_n$ and
it yields a $4\Bbb Z^n$-periodic  cube tiling with no two 
cubes having a
complete facet in common if and only if $\scr S$ forms a 
clique in $G^*_n$.
This gives the Corr\'adi-Szab\'o criterion that a 
Szab\'o-type counterexample
exists in $\Bbb R^n$ if and only if $G^*_n$ contains a 
clique of size $2^n$. 

The graph $G_n$ is the complement of the product graph 
$C_4 \otimes
C_4\otimes\cdots \otimes C_4$ of $n$ copies of the 4-cycle 
$C_4$. It has
maximal clique size equal to the independent set number of 
$C_4\otimes \cdots
\otimes C_4$, which is $\alpha(C_4)^n=2^n$ since $C_4$ is 
a perfect graph and $
\alpha(C_4)=2$. In fact, $G_n$ has an enormous number of 
maximal cliques, and
the problem is whether or not any of them remain a clique 
in $G^*_n$. 

Note also that the graphs $G_n$ and $G^*_n$ have large 
groups of automorphisms.
On both graphs one can relabel the vertices $\bold 
m=(m_1,\dotsc, m_n)$ by
relabeling the $i$\snug th coordinate using the group 
generated by the cyclic
permutations (0123) and the 2-cycle (13), and one also can 
permute coordinates.
This generates a group of $8^n n$! automorphisms.
\vskip1.5pc

\midinsert
\tableS
\toptablecaption{{\smc Table 1.} \rm Clique $\scr T$ in 
$G_3$.}
\format\cmath&\cmath&\cmath\endformat
\ruleh
0&0&0\\
2&0&1\\
1&2&0\\
0&1&2\\
2&0'&3\\
3&2&0'\\
0'&3&2\\2&2&2
\endtableS
\endinsert

\demo{Proof of Theorem \RM A} We give the easier 
$12$-dimensional construction
first. It starts with the set $\scr T$ of vectors given in 
Table 1,
which (ignoring the primes on some zeros) is a clique of 
size 8 in $G_3$. It is
very nearly a clique for $G^*_3$, in that it omits only 
three edges, namely,
201--$20'3$, 120--$320'$, 012--$0'32$. If somehow $0'$ 
were distinct from 0,
while it was still the case that $2-0'=2$, then this would 
be a $2^3$-clique
for $G^*_3$ and would give a counterexample. \enddemo

We use a block substitution construction that in effect 
accomplishes this in $
\Bbb R^{3k}$ for suitable $k$. Assign to each of 
$0,0',1,2,3$ sets $S_0, S'_0,
S_1, S_2, S_3$ of vectors in $\{0,1,2,3\}^k$ having the 
following properties: 
\roster
\item"(i)" Each of $S_0, S'_0, S_1, S_2, S_3$ is a clique 
in $G^*_k$. 
\item"(ii)" No two of the sets $S_0, S'_0, S_1, S_2, S_3$ 
have a common vector. 
\item"(iii)" \<$S_0\cup S_2, S'_0 \cup S_2$, and $S_1\cup 
S_3$ are each a clique
in $G_k$. 
\endroster
Assuming (i), (ii), the last condition (iii) says, e.g., 
for $S_0\cup S_2$ each
element of $S_0$ differs from each element of $S_2$ by 2 
(mod 4) in some
coordinate. 

Call the vectors in the sets $S_i$ blocks. Form the set 
$\scr S$ of all vectors
in $\{0,1,2,3\}^{3k}$ that can be formed by taking any 
vector $(m_1, m_2,
m_3)\in \scr T$ and for each $m_i$ substituting any block 
in the corresponding
$S_{m_i}$, independently for each $i$. 

\rem{Claim} $\scr S$ is a clique in $G^*_{3k}$. 
\endrk 

To prove this, let $\bold v,\bold v'$ be distinct elements 
of $\scr S$
constructed from $\bold m=(m_1,m_2, m_3)$ and $\bold 
m'=(m'_1,m'_2, m'_3)$ in $
\scr T$, respectively. If $\bold m=\bold m'$ then $\bold 
v,\bold v'$ have some
block $\bold w,\bold w'\in S_{m_i}$ where they differ, and 
condition (i) forces
an edge between $\bold v$ and $\bold v'$ in $G^*_{3k}$. If 
$\bold m\neq \bold
m'$ then $\bold m$ and $\bold m'$ differ by 2 in some 
coordinate (here $0'$ and
2 are considered to differ by 2), which carries over to 
$\bold v$ and $\bold
v'$ by condition (iii), and $\bold m$ and $\bold m'$ also 
differ in another
coordinate, where 0 is treated as distinct from $0'$, and 
this carries over to
$\bold v$ and $\bold v'$ by condition (ii), proving the 
claim. 

If one can choose $|S_0|=|S'_0|=a$, $|S_1|=b$, $|S_2|=c$, 
and $|S_3|=d$ with
$a+c=2^k$, $b+d=2^k$, then $|\scr S|=a^3+3abc+3acd+
c^3=2^{3k}$ will be a clique in
$G^*_{3k}$, thus giving a counterexample. 

\midinsert
\tableS
\toptablecaption{{\smc Table 2.} \rm Blocks used in 
constructions.}
\format\cmath&\cmath&\cmath&\cmath&\cmath&\cmath\endformat
\tch\cmath{S_0}&\tch\cmath{S_0'}&\tch 
\cmath{S_2}&\tch\cmath{S_1}
&\tch\cmath{S'_1}&\tch\cmath{S_3}\\
\ruleh
0000&0303&0211&1000&1303&1211\\
0012&1011&1132&1012&2011&2132\\
0213&1113&2303&1213&2113&3303\\
0230&1130&3020&1230&2130&0020\\
0332&1323&&1332&2323&\\
1020&1331&&2020&2331&\\
2100&2211&&3100&3211&\\
2112&3001&&3112&0001&\\
2220&3022&&3220&0022&\\
2301&3103&&3301&0103&\\
2322&3223&&3322&0223&\\
3132&3231&&0132&0231&
\endtableS
\endinsert

We achieve this with $k=4$, $a=b=12$, $c=d=4$, with the 
sets $S_0,S'_0, S_1,
S_2, S_3$ given in Table 2  below. The sets $S_0, S'_0, 
S_2$ were obtained from
a 28-clique for $G^*_5$ given in \cite{1, Table 2}. 
Examining the first
column of this 28-clique, one finds twelve vectors each 
having value
0 and 2 and four having
value 1. Deleting this column and grouping the resulting 
$\Bbb Z^4$-vectors as
$S_0, S'_0, S_2$, the $G^*_5$-clique property guarantees 
that all the
conditions
(i), (ii), (iii) that concern only $S_0, S_0', S_2$, 
automatically hold. Next
we apply a suitable automorphism of $G^*_4$ to $S_0, S_2$ 
to obtain $S_1, S_3$.
For {\it any} automorphism conditions (i) and (iii) will 
automatically hold for
$S_1, S_3$  obtained this way. Thus we need only to find 
an automorphism where
(ii) holds. The automorphism that cyclically permutes the 
labels of the first
coordinate $0\to 1\to 2\to 3\to 0$ gives suitable $S_1, 
S_3$, as listed in
Table 2. 

The conditions (i), (ii), (iii) can be verified directly 
for $S_0, S_0', S_1,
S_2, S_3$ by  hand calculation. Aside from the 
distinctness of all elements,
the automorphism sending $(S_0, S_2)$ to $(S_1, S_3)$ 
means that one need only
check properties for $S_0, S'_0, S_2$. The calculation can 
be further reduced
by observing that there is an automorphism of $G^*_4$ that 
fixes $S_2$ and
sends $S_0$ to $S_0'$. This automorphism cyclically 
permutes the labels of the
first coordinate $0\to 1\to 2\to 3\to 0$ and the last 
coordinate $0\to 3\to
2\to 1\to 0$, and then exchanges these coordinates. Thus 
one need only verify
that $S_0$ and $S_2$ are $G^*_4$-cliques and $S_0\cup S_2$ 
is a $G_4$-clique. 

The 10-dimensional construction is similar in nature and 
is based on the fact
that the set $\widetilde {\scr T}=S_0\cup S_2$ from Table 2 
is a clique of size
$2^4$ in  $G_4$, which is very nearly a $2^4$-clique for 
$G^*_4$. In $G^*_4$ it omits only the four edges 
0213--0211,
3132--1132, 2301--2303, 1020--3020. Now regard $S_2$ as 
being
$$
\matrix
0&2&1'&1\\
1&1'&3&2\\
2&3&0'&3\\
3&0'&2&0
\endmatrix
$$
where we want $0\neq 0'$ and $1\neq 1'$. Assign to $0,0', 
1,1', 2,3$ the sets
of blocks $S_0, S_0', S_1, S'_1, S_2, S_3$ in Table 2, 
where $S_1'$ is
constructed from $S_1$ similarly to $S'_0$ from $S_0$. 
These sets satisfy: 
\roster
\item"(i)" Each of $S_0, S_0', S_1, S'_1, S_2, S_3$ is a 
clique in $G^*_4$. 
\item"(ii)" No two of these sets have a common vector. 
\item"(iii)" \<$S_0\cup S_2, S'_0\cup S_2, S_1\cup S_3$, 
and $S'_1\cup S_3$ are
each a clique in $G_4$. 
\endroster 

Apply the block substitution construction to the second 
and third columns only
on $\widetilde {\scr T}$ to obtain a set $\widetilde {\scr S}$ 
of $2^{10}$
10-vectors. This is a clique in $G^*_{10}$, as required. 
Note that only the
second and third columns need to be expanded in blocks, 
because the primed
elements in $S_2$ above appear only in these columns. \qed

\heading 3. Discussion\endheading

The failure of Keller's Conjecture in high dimensions 
illustrates the general
phenomenon that Euclidean space allows more freedom of 
movement in high
dimensions than in low ones. It is interesting that the 
critical dimension
where Keller's Conjecture first fails, which is at least 
7, is as high as it
is. 

It may be a difficult matter to determine exactly the 
critical dimension.
Exhaustive search for Szab\'o-type counterexamples  
already seems infeasible
for $G^*_7$; the maximum clique problem is a well-known 
NP-complete problem,
which is also computationally hard in practice. The 
authors ruled out the
existence of any $2^7$-clique in $G^*_7$ that is invariant 
under a cyclic
permutation of coordinates by computer search. It is 
conceivable that there
exist Szab\'o-type counterexamples in dimension 7, 8, or 
9, which are all so
structureless that they will be hard to find. In any case 
we have so far found
no variant of the constructions of Theorem A that work in 
these dimensions. 

A natural extension of Keller's conjecture is to determine 
the largest integer
$K_n$ such that every tiling of $\Bbb R^n$ by unit cubes 
contains two cubes
that have a common face of at least dimension $K_n$. For a 
Szab\'o-type tiling,
two cubes having coordinates $(m_1,\dotsc, m_n)$ and 
$(m'_1,\dotsc, m'_n)$ in
$G^*_n$ have a $k$-dimensional face in common if 
$|m_i-m'_i|=0$ or 2 for all
$i$, and exactly $k$ values $|m_i-m'_i|=0$. The 
10-dimensional and
12-dimensional cube tilings $\widetilde{\scr S}$ and $\scr 
S$ constructed in Theorem
A each contain two cubes sharing a common face of 
codimension 2, so they imply
only $K_{10}\le 8$ and $K_{12}\le 10$. We have found a 
different 10-dimensional
cube tiling (using a similar construction) which shows 
that $K_{10}\le 7$. We
also can show that $n-K_n\to \infty$ as $n\to\infty$; 
details will appear elsewhere. 

\heading Acknowledgment\endheading
We thank Victor Klee for bringing  this problem and the 
work of Corr\'adi and
Szab\'o to our attention.

\enddocument